\documentclass[11pt,letterpaper]{article}%
\usepackage{amsfonts}
\usepackage[left=1in,top=1in,right=1in,bottom=1in]{geometry}
\usepackage{setspace}
\usepackage{bbm}
\usepackage{amsmath,amsthm,amssymb}
\usepackage{graphicx}
\usepackage{xcolor}
\usepackage{dsfont}
\usepackage{amsmath}
\usepackage{amssymb}
\usepackage{mathtools}%
\providecommand{\U}[1]{\protect\rule{.1in}{.1in}}
\renewcommand{\Pr}{\mathbb{P}}

\newtheorem{theorem}{Theorem}
\newtheorem{remark}{Remark}

\newtheorem{lemma}{Lemma}
\newtheorem{proposition}{Proposition}

\begin{document}

\author{Paulo Manrique\thanks{Department of Probability and Statistics, CIMAT,
Guanajuato, Mexico, \textit{paulo.manrique@cimat.mx}}
\and Victor P\'{e}rez-Abreu \thanks{Department of Probability and Statistics,
CIMAT, Guanajuato, Mexico, \textit{pabreu@cimat.mx}}
\and Rahul Roy \thanks{Indian Statistical Institute, New Delhi, India,
\textit{rahul@isid.ac.in}}}
\title{\textbf{On the Universality of the Non-Singularity of General Ginibre and
Wigner Random Matrices}}
\date{\today}
\maketitle

\begin{abstract}
We prove the universal asymptotically almost sure non-singularity of general
Ginibre and Wigner ensembles of random matrices when the distribution of the
entries are independent but not necessarily identically distributed and may
depend on the size of the matrix. These models include adjacency matrices of
random graphs and also sparse, generalized, universal and banded random
matrices. We find universal rates of convergence and precise estimates for the
probability of singularity, which depend only on the size of the biggest jump
of the distribution functions governing the entries of the matrix and not on
the range of values of the random entries. Moreover, no moment assumptions are
made about the distributions governing the entries. Our proofs are based on a
concentration function inequality due to Kolmogorov, Rogozin and Kesten, which
allows us to improve universal rates of convergence for the Wigner case when
the distribution of the entries do not depend on the size of the matrix.

\textit{Key terms: }Adjacency matrix of random graphs, banded random matrix,
decoupling, concentration function, generalized Wigner ensemble,
Littlewood--Offord inequality, Kolmogorov--Rogozin inequality, nondegenerate
distribution, sparse random matrix.

\end{abstract}

\section{Introduction and main results}

Let $A_{n}=($ $\xi_{ij}^{(n)})$ be an $n\times n$ random matrix where each
entry $\xi_{ij}^{(n)}$ follows a distribution $F_{ij}^{(n)}$, $1\leq i,j\leq
n$. The study of the non-singularity of such matrices has mainly been
considered when $F_{ij}^{(n)}\equiv F$ and for two ensembles of random
matrices, the Ginibre and Wigner. We will use the following terminology: An
$n\times n$ random matrix $G_{n}=\left(  \xi_{ij}\right)  _{1\leq i,j\leq n}$
is called a Ginibre matrix if $\xi_{ij}$, $i,j=1,\ldots,n$ are independent
random variables, and an $n\times n$ random symmetric matrix $W_{n}=\left(
\xi_{ij}\right)  _{1\leq i,j\leq n}$ is called Wigner matrix if $\xi_{ij}%
=\xi_{ji}$, $i,j=1,\ldots,n$ and $\xi_{ij}$, $1\leq i\leq j\leq n$ are
independent random variables. We will not assume that the distributions of the
entries have moments.

The singularity of these matrices is trivial if the distributions of $\xi
_{ij}$ are degenerate. The non-singularity is also straightforward if the
entries have continuous distributions. The interesting situation occurs when
some of the entries have distributions with jumps. The singularity of such
matrices is a highly non-trivial problem.

The study of the non-singularity of Ginibre matrices goes back to the
pioneering work by Koml\'{o}s. In \cite{Komlos1967} he considers Ginibre
random matrices $GB(n,1/2)$, whose entries are i.i.d. Bernoulli random
variables, taking the values $0$ or $1$ with probability $1/2$ each. Using a
very clever `growing rank analysis' together with the Littlewood--Offord
inequality, Koml\'{o}s proved that\newline$\mathbb{P}\left\{
\mbox{rank}(GB(n,1/2))<n\right\}  =\mbox{o}(1)$ as $n$ $\rightarrow\infty$.
Bollob\'{a}s \cite{Bollobas1985} presents the concept of `strong rank' and
together with the Littlewood--Offord inequality obtains an unpublished result
due to Koml\'{o}s, \textit{viz.} $\mathbb{P}\left\{
\mbox{rank}(GB(n,1/2))<n\right\}  =\mbox{O}(n^{-1/2})$ as $n$ $\rightarrow
\infty$. Koml\'{o}s \cite{Komlos1968} was also the first to consider the
singularity of Ginibre matrices whose entries are i.i.d. random variables with
a \textit{common} arbitrary non-degenerate distribution, proving that the
probability that such an $n\times n$ matrix is singular has order
$\mbox{o}(1)$ as $n$ $\rightarrow\infty$. This result was improved by Kahn,
Koml\'{o}s and Szemer\'{e}di \cite{KahnKomlosSzemeredi1995} in the case of
Ginibre matrices whose entries are i.i.d. taking values $-1$ or $1$ with
probability $1/2$ each, showing that the probability of singularity is bounded
above by $\theta^{n}$ for $\theta=.999$. The value of $\theta$ has been
improved by Tao and Vu \cite{TaoVu2006}, \cite{TaoVu2007} to $\theta=3/4+o(1)$
and by Bourgain, Vu and Wood \cite{BoVuWo10} to $\theta=1/\sqrt{2}+o(1).$
Slinko \cite{Si01} considered Ginibre random matrices whose entries have the
same uniform distribution taking values in a finite set, proving also that the
probability of singularity is $\mbox{O}(n^{-1/2})$ as $n$ $\rightarrow\infty$.

The aim of this paper is to understand the asymptotic non-singularity of more
general Ginibre and Wigner ensembles. We are interested in finding
universality results with respect to general distributions of the entries and
also when these distributions depend on the size of the matrix.

As a first step in this direction, the results in \cite{Bollobas1985},
\cite{Si01} were generalized by Bruneau and Germinet \cite{BrGe09} to Ginibre
random matrices whose entries follow different independent non-degenerate
distributions $F_{ij}$ which do not change with the size of the matrix. Their
result gives a universal rate of convergence of $n^{-1/2}$ as follows:

\begin{proposition}
\label{BrueauGerminet}\ (Bruneau and Germinet, 2009). Let $G_{n}$ be an
$n\times n$ Ginibre matrix with independent entries $\xi_{ij}$ satisfying the
following property $H$: there exists $\rho\in(0,1/2)$ such that for any
$i,j=1,\ldots,n$, $\mathbb{P}\left\{  \xi_{ij}>x_{ij}^{+}\right\}  >\rho$ and
$\mathbb{P}\left\{  \xi_{ij}<x_{ij}^{-}\right\}  >\rho$ \ for some real
numbers $x_{ij}^{-}<x_{ij}^{+}$, then
\begin{equation}
\mathbb{P}\left\{  \mbox{rank}(G_{n})<n\right\}  \leq C/\sqrt{\beta_{\rho
}(1-\rho)n}, \label{BrueGermresult}%
\end{equation}
where the constant $C$ is universal (coming from the Littlewood--Offord
inequality) and $\beta_{\rho}$ is an implicit constant $0<\beta_{\rho}<1$
which goes to zero as $\rho\rightarrow1$.
\end{proposition}

\begin{remark}
a) The above proposition is proved in \cite{BrGe09} using ideas of strong rank
of \cite{Bollobas1985}, together with a Bernoulli representation theorem for
the distribution of a random variable, and the Littlewood--Offord inequality.

b) We point out that it is possible to express (\ref{BrueGermresult}) in terms
of the size of the biggest jump of the distribution functions governing the
entries. Indeed, this follows using a strong rank analysis and the
Kolmogorov--Rogozin concentration inequality. This inequality, stated in
Section 2, will be used repeatedly in this work. Returning to
(\ref{BrueauGerminet}), taking $\kappa=\max_{1\leq i,j\leq n}\sup
_{x\in\mathbb{R}}\Pr\{\xi_{i,j}=x\}$, the size of the biggest jump of $F_{ij}%
$, $i,j=1,\ldots,n$, for $0\geq\kappa<1$ we have
\begin{equation}
\mathbb{P}\left\{  \mbox{rank}(G_{n})<n\right\}  \leq\frac{C_{1}}{\sqrt
{\beta_{\kappa}^{\prime}(1-\kappa)n}}, \label{OurResultGinibreStrongRank}%
\end{equation}
where the constant $C_{1}$ is a universal constant coming from the Kolmogorov--Rogozininequality.

c) We observe that the constants $\beta_{\rho}$ and $\beta_{\kappa}^{\prime}$
are not universal: they might depend on the distributions $F_{ij}$.

d) These two results highlight the fact that the non-singularity of Ginibre
matrix depends only on $\rho$ or, equivalently, the size of the biggest jump
$\kappa$. In other words, the universal property of a random matrix being
non-singular depends neither on the range of values taken by the entries nor
on other properties of their distribution except the size of the biggest jump.
\end{remark}

As for Wigner random matrices, the study of their singularity was initiated by
Costello, Tao and Vu \cite{CostelloTaoVu2006} inspired by the work of
Koml\'{o}s \cite{Komlos1967}.

\begin{proposition}
\label{CostTaiVu} (Costello et al. 2006). Let $W_{n}=(\xi_{ij})$ be an
$n\times n$ Wigner matrix whose upper diagonal entries $\xi_{ij}$ are
independent random variables with common Bernoulli distribution on $\left\{
0,1\right\}  $ with parameter $1/2$. Then, as $n\rightarrow\infty$,
\[
\mathbb{P}\left\{  \mbox{rank}(W_{n}))<n\right\}  =\mbox{O}(n^{-1/8+\alpha}),
\]
for any positive constant $\alpha$, the implicit constant in O$(\cdot)$
depending on $\alpha$.
\end{proposition}

\begin{remark}
a) The proof of the above proposition in \cite{CostelloTaoVu2006} required
developing a quadratic Littlewood--Offord inequality. A possible
generalization to distributions other than Bernoulli was also indicated in
\cite{CostelloTaoVu2006}.

b) Theorem \ref{Univer}.b below gives a better universal rate of convergence
$n^{-1/4+\alpha}$, for any Wigner random matrix $W_{n}=(\xi_{ij})$ with
independent entries, which need not be identical. While the off-diagonal
entries need to be non-degenerate, the diagonal entries could be degenerate.
\end{remark}

More recently, Wigner matrices have been studied when the entries satisfy some
restrictions. Nguyen \cite{Nguyen2013} considered a Wigner matrix $W_{n} $
with entries taking the values $-1$ or $1$ with probability $1/2$ each,
subject to the condition that each row has exactly $n/2$ entries which are
zero. He showed that the probability of $W_{n}$ being singular is
$\mbox{O}(n^{-C}),$ for any positive constant $C$, the implicit constant in
O$(\cdot) $ depending on $C$. Recently, Vershynin \cite{Ve12} has considered
the case of a Wigner matrix $W_{n}$ whose entries satisfy the following
property: the above-diagonal entries are independent and identically
distributed with zero mean, unit variance and subgaussian, while the diagonal
entries satisfy $\xi_{ii}\leq K\sqrt{n}$ for some $K$. He showed that the
probability of $W_{n}$ being singular is bounded above by $2\exp(-n^{c}),$
where $c$ depends only on the subgaussian distribution and on $K$.

One of the goals of this paper is to study the non-singularity of Ginibre and
Wigner matrices when the distributions of \ the entries $F_{ij}^{(n)}$ depend
on the size of the matrix. This kind of random matrix appears in the study of
random graphs \cite{CostelloVu2006}, circular law \cite{BoCh12}, sparse
matrices \cite{CoVu10}, \cite{Er2010} and some other models that have recently
been extensively considered, such as the so-called generalized, universal and
banded Wigner ensembles \cite{ErYaYi12}, \cite{Sp11}, among other papers. See
also the non-i.i.d. Wigner case in, for example, \cite[pp 26]%
{BaiSilverstein2006}.

One difficulty that arises in this situation is to find adequate asymptotic
estimates for the probability of the singularity's being zero, such that the
constants involved in the rate of convergence do not depend on the
distributions of the entries. We overcome this difficulty using a universal
concentration inequality due to Kesten \cite{Kesten1972}, which we express in
terms of the size of the jumps of the distribution functions.

\paragraph{1.1 Main results}

We now consider Ginibre and Wigner matrix ensembles $G_{n}^{(n)}=\left(
\xi_{ij}^{(n)}\right)  _{1\leq i,j\leq n}$, and $W_{n}^{(n)}=\left(  \xi
_{ij}^{(n)}\right)  _{1\leq i,j\leq n} $, where the distribution function
$F_{ij}^{(n)}$ governing $\xi_{ij}^{(n)}$ is allowed to change with the size
of the matrix.

One of our main conclusions is as to the non-singularity of the above Ginibre
and Wigner ensembles. More specifically, given a collection of non-degenerate
distribution functions $\{F_{ij}^{(n)}:\;i,j\geq1,\;n\geq1\}$ and a
subsequence $\{m_{n}:n\geq1\}$, we study the singularity of the $m_{n}\times
m_{n}$ matrix with independent entries $\xi_{kl}^{(n)}$ governed by the
distribution function $F_{kl}^{(n)}$ for every $1\leq k,l\leq m_{n}$. Let us
denote by $\kappa_{n}$ the biggest jump of the distribution functions
$F_{ij}^{(n)}$, $1\leq i,j\leq m_{n}$, i.e., if $\kappa_{i,j}=\sup
_{x\in\mathbb{R}}\Pr\{\xi_{i,j}^{n}=x\}$, then
\begin{equation}
\kappa_{n}=\max_{1\leq i,j\leq n}\{\kappa_{i,j}\}. \label{BiggestJump}%
\end{equation}
We give a sufficient condition for $m_{n}=n$ in terms of the sequence of
biggest jumps $\left(  \kappa_{n}\right)  _{n\geq1}$.

\begin{theorem}
\label{Univer} (Universality of the non-singularity of Ginibre and Wigner
ensembles) With the notation as above, let $G_{r}^{(n)}$ and $W_{r}^{(n)}$ be
the $r\times r$ Ginibre and Wigner matrices respectively, each with entries
$\xi_{i,j}^{(n)}$, $1\leq i\leq j\leq r$. Assume that $\kappa_{n}<\kappa
\in[0,1)$ for all $n$

a) As $n\to\infty$
\begin{equation}
\mathbb{P}\left\{  \mbox{rank}(G_{n}^{(n)})<n\right\}  =\mbox{O}\left(
n^{-1/2}\right)  \label{eq:Univer11}%
\end{equation}
where the implicit constant in O$(\cdot)$ depends on $\kappa$.

b) For any $\varepsilon\in(0,1)$,
\begin{equation}
\mathbb{P}\left\{  \mbox{rank}(W_{n}^{(n)})<n\right\}  =\mbox{O}\left(
n^{-(1-\varepsilon)/4}\right)  , \label{eq:Univer21}%
\end{equation}
where the implicit constant in O$(\cdot)$ depends on $\varepsilon$ and
$\kappa$.
\end{theorem}

\begin{remark}
a) The proof of Lemma \ref{CorGin01} in Section 3, where Theorem
\ref{Univer}.a is given, highlights the fact that the probability that a
Ginibre matrix has small rank is small. This is due only to the independence
of entries.

b) The bound $n^{-1/4+\alpha}$ in (\ref{eq:Univer21}) improves the rate
$n^{-1/8+\alpha}$ in Theorem \ref{CostTaiVu} of Costello, Tao and Vu
\cite{CostelloTaoVu2006}.
\end{remark}

We now turn to Theorem \ref{Univer}. A natural question is to understand what
happens when $\kappa_{n}\to1$.

\begin{proposition}
\label{ProGin01} For any sequence $\{\kappa_{n}\in[0,1] : n\geq1\}$ there is a
sequence $\{G_{m_{n}}=(\xi_{i,j})_{1\leq i,j\leq m_{n}}\}$ such that:

\begin{itemize}

\item $G_{m_{n}}$ is a $m_{n}\times m_{n}$ Ginibre matrix

\item $\xi_{i,j}$, $1\leq i,j\leq m_{n}$, have the same distribution
$F_{m_{n}}$

\item $\kappa_{n}$ is the maximum jump of $F_{m_{n}}$

\item $\Pr\{ G_{m_{n}} \mbox{ has full rank }\}\to1\;\;\;\;\;n\to\infty$
\end{itemize}
\end{proposition}

In the following examples we can see that if $\kappa_{n}\rightarrow1$ at some
appropiate rate, the probability of a singularity can behave differently.

We write $GB(n,p)$ ($WB(n,p)$) for a $n\times n$ Ginibre (Wigner) matrix whose
entries obey a Bernoulli distribution on $\{0,1\}$ with parameter $p$.

Let $ZGB_{n}$ ($ZWB_{n}$) be the event that the first row of $GB(n,1/n)$,
($WB(n,1/n)$) contains only zeros. Then
\[
\mathbb{P}\left\{  ZGB_{n}\right\}  =\left(  1-\frac{1}{n}\right)
^{n},\;\;\;\;\;\mathbb{P}\left\{  ZWB_{n}\right\}  =\left(  1-\frac{1}%
{n}\right)  ^{n},
\]
and hence
\[
e^{-1}\leq\lim_{n\rightarrow\infty}\mathbb{P}\left\{  \mbox{rank}\left(
GB\left(  n,1/n\right)  \right)  <n\right\}  ,
\]%
\[
e^{-1}\leq\lim_{n\rightarrow\infty}\mathbb{P}\left\{  \mbox{rank}\left(
WB\left(  n,1/n\right)  \right)  <n\right\}  .
\]

However, if $\alpha\in(0,1)$, then there is a constant $C_{\alpha}>0$
%\begin{equation}
%\mathbb{P}\left\{  \mbox{rank}\left(  GB\left(  n,n^{\alpha}/n\right)
%\right)  <n\right\}  \leq C  n^{-\alpha/4}  , \label{Sd3}%
%\end{equation}%
\begin{equation}
\mathbb{P}\left\{  \mbox{rank}\left(  WB\left(  n,n^{\alpha}/n\right)
\right)  <n\right\}  \leq n^{-C_{\alpha}}. \label{Sd4}%
\end{equation}
In the Ginibre case it is not clear what happens when $\kappa=n^{\alpha}/n$,
but if $\gamma\in(0,1)$, then
\begin{equation}
\label{Sd5}\Pr\{\mbox{rank}\left(  GB\left(  n,n^{\alpha}/n\right)  \right)
>\gamma n \} \to1\;\;\mbox{as}\;\;n\to\infty.
\end{equation}

Furthermore, as an application of the Wigner case, we obtain an estimation of
the probability that the adjacency matrix of a sparse random graph (not
necessarily an Erd\"{o}s--R\'{e}nyi graph) is non-singular. Costello and Vu
\cite{CostelloVu2006} have analyzed the adjacency matrices of sparse
Erd\"{o}s--R\'{e}nyi graphs where each entry is equal to 1 with the same
probability $p(n)$, which tends to 0 as $n$ goes to infinity (see also
Costello and Vu \cite{CoVu10}, where a generalization of \cite{CostelloVu2006}
is considered in which each entry takes the value $c\in\mathbb{C}$ with
probability $p$ and zero with probability $1-p,$ and the diagonal entries are
possibly non-zero). It is proved in \cite{CostelloVu2006} that when
$c\ln(n)/n\leq p(n)\leq1/2$, $c>1/2$, then with probability $1-$O$((\ln
\ln(n))^{-1/4})$, the rank of the adjacency matrix equals the number of
non-isolated vertices. Now we consider the following model extension of
Erd\"{o}s--R\'{e}nyi graphs, where vertices $i$ and $j$ are linked with a
probability that depends on $i$ and $j$ and the number of vertices.
Furthermore, the rate of convergence is an improvement of the one given in
\cite{CostelloVu2006} for $c\ln n/n^{\beta}\leq p(n)\leq1/2$ with $c>0$ and
$\beta\in(0,1)$. From the proof of Theorem \ref{Univer}.b in Section 4, if
$\kappa_{n}=1-p(n)$, we have as $n\to\infty$
\[
\frac{\kappa_{n}^{\frac{3}{8}n-\frac{1}{2}n^{1-\varepsilon}}} {\kappa
_{n}(1-\kappa_{n})} \leq\left(  \frac{\kappa_{n}^{2}}{n^{1-\varepsilon
}(1-\kappa_{n})}\right)  ^{1/4}\leq\left(  \frac{(1-c(\ln n/n^{\beta}))^{2}%
}{n^{1-\varepsilon-\beta} \ln n}\right)  ^{1/4}\to0
\]
if $\varepsilon+\beta<1$.

\begin{proposition}
\label{RandGrap} Let $\{p_{ij}\in(0,1): i,j=1,2,\ldots\}$ be a double sequence
of positive numbers with $p_{n}^{\ast}=\min_{1\leq i\leq j\leq n}\{p_{ij}%
\}\in[c\ln n/n^{\beta},1/2]$, $c>0$, and $\varepsilon+\beta<1$, $\varepsilon
,\beta\in(0,1)$. Then there is a random graph with $n$ vertices such that the
vertex $i$ is linked with the vertex $j$ with probability $p_{ij}$, $1\leq
i<j\leq n$, and if $A_{n}$ is the adjacency matrix, we have as $n\to\infty$
\begin{equation}
\Pr\left\{  \mbox{rank}(A_{n})<n\right\}  \leq Cn^{-(1-\varepsilon-\beta)/4},
\label{ProEq02}%
\end{equation}
for some constant $C>0$.
\end{proposition}

\begin{remark}
a) In many applications of random matrices one considers ensembles of the form
$G_{n}^{(n)}=a_{n}^{-1}G_{n}$ and $W_{n}^{(n)}=a_{n}^{-1}W_{n}$ where
$a_{n}\rightarrow\infty$ as $n\rightarrow\infty$ and the non-degenerate
distributions of the entries of $\ G_{n}$ and $W_{n}$ do not depend on the
matrix size, $n.$ In this case $\kappa_{n}=\kappa<1$ for all $n\geq1$ if the
distribution is not degenerate. However the ensembles $G_{n}^{(n)}$ and
$W_{n}^{(n)}$ are asymptotically almost surely non-singular. In fact, this
holds for any sequence $a_{n}\rightarrow\infty$ and the rate of convergence to
zero of the probability of singularity is not affected by the rate of
convergence of $a_{n}$ if the distributions of the entries have discrete support.

b) The case $a_{n}=1/\sqrt{n}$ is the setup of those problems of random
matrices appearing in the study of asymptotic spectral distributions
\cite{AndGuitZei2010}, \cite{BaiSilverstein2006}, geometric functional
analysis \cite{Sa90}, \cite{Ru08}, and restricted isometries \cite{RuVe10},
among others.

c) Finally, the results in the Ginibre case have a straightforward extension
to non-square $n\times m$ random matrices whose entries are independent random
variables and have distributions with jumps.
\end{remark}

\section{Preliminaries on Concentration Inequalities}

\label{Section:Preliminaries}

In this section we present the Kolmogorov--Rogozin concentration inequalities
that we use for the proofs of our main results on non-singularity. We express
these inequalities in terms of the size of the biggest jump of the
non-degenerate distribution functions.

The L\'{e}vy concentration function $Q(\xi;\lambda)$ of a random variable
$\xi$ is defined by
\[
Q(\xi;\lambda)=\sup_{x\in\mathbb{R}}\Pr\left\{  \xi\in\left[  x,x+\lambda
\right]  \right\}  ,\;\;\;\lambda>0.
\]
Let $\xi_{1},\xi_{2},\ldots$ be independent random variables and $S_{n}%
=\sum_{i=1}^{n}\xi_{i}$. An expression that relates the concentration function
of $S_{n}$ to the concentration functions of the summands $\xi_{i}$ was given
by Kolmogorov--Rogozin; see \cite{Kesten1969}.

\begin{lemma}
[The Kolmogorov--Rogozin Inequality]\label{LemKolRog} There exists a universal
constant $C$ such that for any independent random variables $\xi_{1}%
,\ldots,\xi_{n}$ and any real numbers $0<\lambda_{1},\ldots,\lambda_{n}\leq
L$, one has
\[
Q(S_{n};L)\leq CL\left\{  \sum_{i=1}^{n}\lambda_{i}^{2}\left[  1-Q(\xi
_{i};\lambda_{i})\right]  \right\}  ^{-1/2}.
\]

\end{lemma}

Kesten \cite{Kesten1972} obtained the following refinement of the above inequality.

\begin{lemma}
\label{LemKes} For the constant $C$ of the Kolmogorov--Rogozin inequality and
any independent random variables $\xi_{1},\ldots,\xi_{n}$, and real numbers
$0<\lambda_{1},\ldots,\lambda_{n}\leq2L$, one has
\[
Q(S_{n};L)\leq4\cdot2^{1/2}(1+9C)L\frac{\sum_{i=1}^{n}\lambda_{i}^{2}\left[
1-Q(\xi_{i};\lambda_{i})\right]  Q(\xi_{i},L)}{\left\{  \sum_{i=1}^{n}%
\lambda_{i}^{2}\left[  1-Q(\xi_{i};\lambda_{i})\right]  \right\}  ^{3/2}}.
\]

\end{lemma}

For the study of the non-singularity of random matrices, one has to find an
estimate of the probability that a polynomial of independent random variables
equals a real number. In the case of Ginibre and\ Wigner matrices, the
polynomials are of degree one and two, respectively. Our first goal is to
write the Kesten inequality in terms of the size of the biggest jump and then
obtain the corresponding linear and quadratic concentration inequalities.

We first discuss the relation between the size of the biggest jump of a
non-degenerate distribution $F$ and its corresponding L\'{e}vy concentration
function. Let $D_{F}$ be the set of discontinuities of $F$ and $\kappa$ its
biggest jump, i.e., $\kappa=\sup_{x\in\mathbb{R}}\Pr\left\{  \xi=x\right\}  $,
where $\xi$ has the distribution function $F$.

We note the following:

\begin{enumerate}
\item There exists $x_{\kappa}\in\mathbb{R}$ such that $\Pr\left\{
\xi=x_{\kappa}\right\}  =\kappa$.

\item Let $p_{i}=\Pr\left\{  \xi=x_{i}\right\}  $, $i\in\mathbb{N}$, then
$\sum_{i\geq1}p_{i}\leq1$, i.e., for all $\varepsilon>0$ there exists
$N(\varepsilon)\in\mathbb{N}$ such that $\sum_{i\geq n}p_{i}\leq\varepsilon$
for all $n\geq N(\varepsilon)$.

\item If $F$ is a discrete distribution $\left(  \sum_{i\in\mathbb{N}}%
p_{i}=1\right)  $ and $x_{\kappa}$ is not an accumulation point of $D_{F},$
there exists $\delta_{1}>0$ with
\[
\sup_{x\in\mathbb{R}}\Pr\left\{  \xi\in\left[  x,x+\delta_{1}\right]
\right\}  =\kappa.
\]
Otherwise, if $F$ is not discrete or $x_{\kappa}$ is an accumulation point of
$D_{F}$, there exists some $\Delta>0$, which may be taken as small as desired,
such that, for $\Delta$ fixed, there is a $\delta_{2}>0$ with
\[
\sup_{x\in\mathbb{R}}\Pr\left\{  \xi\in\left[  x,x+\delta_{2}\right]
\right\}  =\kappa+\Delta<1.
\]
We define $\kappa_{\Delta}$, for $\Delta\in(0,1)$ fixed, by $\kappa_{\Delta
}:=\kappa$ if $F$ is discrete and $x_{\kappa}$ is not an accumulation point of
$D_{F}$, and otherwise, $\kappa_{\Delta}:=\kappa+\Delta$. So, we have that
there is a $\delta>0$ such that
\begin{equation}
\sup_{x\in\mathbb{R}}\Pr\left\{  \xi\in\left[  x,x+\delta\right]  \right\}
=\kappa_{\Delta}.\;\;\; \label{supDel}%
\end{equation}

\item We fix $\Delta\in(0,1)$ and $\delta>0$ satisfying (\ref{supDel}). If
$a\in\mathbb{R}$ with $|a|\geq1$, then
\[
\sup_{x\in\mathbb{R}}\Pr\left\{  a\xi\in\left[  x,x+\delta\right]  \right\}
\leq\kappa_{\Delta}.
\]
Indeed, if $\sup_{x\in\mathbb{R}}\Pr\left\{  a\xi\in\left[  x,x+\delta\right]
\right\}  >\kappa_{\Delta}$, then there exists some $x^{\ast}\in\mathbb{R}$
such that \newline$\Pr\left\{  a\xi\in\lbrack x^{\ast},x^{\ast}+\delta
]\right\}  >\kappa_{\Delta}$, but
\[
\delta\geq|a\xi-x^{\ast}|=|a|\left\vert \xi-\frac{x^{\ast}}{a}\right\vert
\geq\left\vert \xi-\frac{x^{\ast}}{a}\right\vert ,
\]
which is a contradiction of the definition of $\kappa_{\Delta}$. So, we have
that
\[
Q(\beta\xi,\delta)\leq\kappa_{\Delta}\mbox{ for } |\beta|\geq1.
\]

\end{enumerate}

Now let $\xi_{1},\ldots,\xi_{n}$ be independent random variables with
distribution functions $F_{1},\ldots,F_{n}$, respectively. For each $\xi_{i}$,
we consider $\kappa(i),\kappa_{\Delta}(i)<1$ defined as above. We first prove
the following concentration inequality in terms of the biggest jumps of the
distribution functions.

\begin{lemma}
[The Linear Concentration Inequality]\label{LemLinCon} Let $\xi_{1},\ldots
,\xi_{n}$ be independent random variables with non-degenerate distributions
$F_{1},\ldots,F_{n}$, respectively, and let $\alpha_{1},\ldots,\alpha_{n}$ be
real numbers with $\alpha_{i}\neq0$, $i=1,\ldots,n$. Then
\[
\sup_{x\in\mathbb{R}}\Pr\left\{  \sum_{i=1}^{n}\alpha_{i}\xi_{i}=x\right\}
=\mbox{O}\left(  \frac{\sum_{i=1}^{n}(1-\kappa(i))\kappa_{\Delta}(i)}{\left\{
\sum_{i=1}^{n}\left[  1-\kappa_{\Delta}(i)\right]  \right\}  ^{3/2}}\right)
,
\]
where the implicit constant in O($\cdot$) does not depend on $F_{i}$,
$i=1,\ldots,n$.
\end{lemma}

\textbf{Proof.} Let $a=\min_{1\leq i\leq n}\left\{  |\alpha_{i}|\right\}  $
and $\delta=\min_{1\leq i\leq n}\left\{  \delta_{i}\right\}  $, where
$\delta_{i}>0$ satisfies $\kappa_{\Delta}(i)= Q(\xi_{i},\delta_{i})$,
$i=1,\ldots,n$. We have for $x\in\mathbb{R}$
\[
\Pr\left\{  \sum_{i=1}^{n}\alpha_{i}\xi_{i}=x\right\}  =\Pr\left\{  \sum
_{i=1}^{n}\frac{\alpha_{i}}{a}\xi_{i}=\frac{x}{a}\right\}  =\Pr\left\{
\sum_{i=1}^{n}\alpha_{i}^{\prime}\xi_{i}=x^{\prime}\right\}  ,
\]
where $\alpha_{i}/a=\alpha_{i}^{\prime}$ and $x/a=x^{\prime}$. Now,
\begin{align*}
\Pr\left\{  \sum_{i=1}^{n}\alpha_{i}^{\prime}\xi_{i}=x^{\prime}\right\}   &
\leq\sup_{y\in\mathbb{R}}\Pr\left\{  \sum_{i=1}^{n}\alpha_{i}^{\prime}\xi
_{i}\in\lbrack y,y+\delta]\right\} \\
&  \leq4\cdot2^{1/2}(1+9C)\frac{\sum_{i=1}^{n}(1-\kappa(i))\kappa_{\Delta}%
(i)}{\left\{  \sum_{i=1}^{n}\left[  1-\kappa_{\Delta}(i)\right]  \right\}
^{3/2}},
\end{align*}
the last expression following from Lemma~\ref{LemKes}. \hfill$\blacksquare$

\begin{remark}
\label{RemLinConIne} a) If $\kappa_{\Delta}(i)<\kappa<1$ for all $i$,
\[
\sup_{x\in\mathbb{R}}\Pr\left\{  \sum_{i=1}^{n}\alpha_{i}\xi_{i}=x\right\}
=\mbox{O}\left(  \frac{\kappa}{\sqrt{(1-\kappa)^{3}n}}\right)  .
\]

b) Lemma \ref{LemLinCon} holds when $r$ of the random variables $\xi
_{1},\ldots,\xi_{n}$ are degenerate for some $1\leq r<n$; in this situation
$n$ is replaced by $n-r$. The contribution to the bound of the concentration
inequality is provided only by the non-degenerate random variables.
\end{remark}

In order to prove the so-called Quadratic Concentration Inequality, we recall
the decoupling argument.

\begin{lemma}
[Decoupling]\label{LemDec} Let $X\in\mathbb{R}^{m_{1}}$ and $Y\in
\mathbb{R}^{m_{2}}$ be independent random variables, with $m_{1}+m_{2}=n$, and
let $\varphi:\mathbb{R}^{n}\rightarrow\mathbb{R}$ be a Borel function. Let
$X^{\prime}$ be a variable independent of $X$ and $Y$, but having the same
distribution as $X$. For any interval $I$ of $\mathbb{R}$, we have
\[
\Pr^{2}\left\{  \varphi(X,Y)\in I\right\}  \leq\Pr\left\{  \varphi(X,Y)\in
I,\varphi(X^{\prime},Y)\in I\right\}  .
\]

\end{lemma}

A quadratic Littlewood--Offord inequality for independent $\{0,1\}$-Bernoulli
random variables with probability $1/2$ was proved in \cite{CostelloTaoVu2006}%
. The result below is for independent random variables not necessarily
identically distributed and without any assumption on their moments.

\begin{lemma}
[The Quadratic Concentration Inequality]\label{LemQuaCon} Let $\xi_{1}%
,\ldots,\xi_{n}$ be independent random variables with non-degenerate
distributions $F_{1},\ldots,F_{n}$, respectively, and let $(c_{ij})_{1\leq
i,j\leq n}$ be a symmetric $n\times n$ array of constants. Suppose
$S_{1}\sqcup S_{2}$ is a partition of $\{1,2,\ldots,n\}$ such that for each
$j\in S_{2}$, the set $N_{j}:=\{i\in S_{1} : c_{ij}\neq0\}$ is non-empty. Let
\[
\varphi=\varphi\left\{  \xi_{1},\ldots,\xi_{n}\right\}  =\sum_{1\leq i,j\leq
n}c_{ij}\xi_{i}\xi_{j}%
\]
be the quadratic form whose coefficients are $c_{ij}$. Then
\[
\Pr\left\{  \varphi=x\right\}  = \mbox{O} \left(  \left[  \frac{1}{|S_{2}%
|}\sum_{j\in S_{2}}\left(  \frac{\sum_{i\in N_{j}}(1-\overline{\kappa
}(i))\overline{\kappa}_{\Delta}(i)}{\left\{  \sum_{i\in N_{j}}\left[
1-\overline{\kappa}_{\Delta}(i)\right]  \right\}  ^{3/2}}\right)  +
\sup_{D\subset S_{2},|D|\geq|S_{2}|/2} \frac{\sum_{j\in D}(1-\kappa
(j))\kappa_{\Delta}(j)}{\left\{  \sum_{j\in D}\left[  1-\kappa_{\Delta
}(j)\right]  \right\}  ^{3/2}}\right]  ^{1/2} \right)
\]
where, for $\xi_{i}^{\prime}$ an independent copy of $\xi_{i}$, $\overline
{\kappa}(i)$ and $\overline{\kappa}_{\Delta}(i)$ are the jumps associated with
$\xi_{i}-\xi_{i}^{\prime}$ and $\kappa(j)$ and $\kappa_{\Delta}(j)$ are the
jumps associated with $\xi_{j}$. The implicit constant in O$(\cdot)$ does not
depend on $F_{i}$, $i=1,\ldots,n$.
\end{lemma}

\textbf{Proof.} Let $\delta=\min_{1\leq i\leq n}\left\{  \delta_{i}\right\}  $
where $\delta_{i}>0$ satisfies $\kappa_{\Delta}(i)=Q(\xi_{i},\delta_{i})$,
$i=1,\ldots,n$. If $x\in\mathbb{R}$, we have
\[
\Pr\left\{  \varphi=x\right\}  \leq\Pr\left\{  \varphi\in\lbrack
x,x+\delta/2]\right\}  .
\]
Write $I=[x,x+\delta/2]$, $X=(\xi_{i}:i\in S_{1})$, $Y=(\xi_{i}:i\in S_{2})$
and $X^{\prime}=(\xi_{i}^{\prime}:i\in S_{1})$, with $X^{\prime}$ independent
of $X$ and $Y$, but having the same distribution as $X$. By Lemma~\ref{LemDec}%
,
\begin{align*}
{\Pr}^{2}\left\{  \varphi(X,Y)\in I\right\}   &  \leq\Pr\left\{
\varphi(X,Y)\in I,\varphi(X^{\prime},Y)\in I\right\} \\
&  \leq\Pr\left\{  \varphi(X,Y)-\varphi(X^{\prime},Y)\in\lbrack-\delta
/2,\delta/2]\right\}  .
\end{align*}
We can rewrite $\varphi(X,Y)-\varphi(X^{\prime},Y)$ as
\begin{align*}
\varphi(X,Y)-\varphi(X^{\prime},Y)  &  = g(X,X^{\prime}) +2\sum_{j\in S_{2}%
}\xi_{j}\left(  \sum_{i\in S_{1}}c_{ij}\left(  \xi_{i}-\xi_{i}^{\prime
}\right)  \right) \\
&  = g(X,X^{\prime})+2\sum_{j\in S_{2}}\xi_{j}\eta_{j},
\end{align*}
where $g(X,X^{\prime})=\sum_{i,j\in S_{1}}c_{ij}(\xi_{i}\xi_{j}-\xi
_{i}^{\prime}\xi_{j}^{\prime})$ and $\eta_{j}=\sum_{i\in S_{1}}c_{ij}\left(
\xi_{i}-\xi_{i}^{\prime}\right)  $.

Let $\zeta$ be the number of $\eta_{j}$ which are equal to zero. If
$J=[-\delta/2,\delta/2]$, we have
\begin{align*}
\Pr\left\{  \varphi(X,Y)-\varphi(X^{\prime},Y)\in J\right\}   &  \leq
\Pr\left\{  \varphi(X,Y)-\varphi(X^{\prime},Y)\in J,\zeta\leq\frac{|S_{2}|}%
{2}\right\} \\
&  +\Pr\left\{  \zeta>\frac{|S_{2}|}{2}\right\}  .
\end{align*}
Since $\zeta=\sum_{j\in S_{2}}\mathbf{1}_{\left\{  \eta_{j}=0\right\}  }$,
using Lemma \ref{LemLinCon}, we have
\begin{align*}
\mathbb{E}\left(  \zeta\right)  =\sum_{j\in S_{2}}\Pr\left\{  \eta
_{j}=0\right\}   &  =\sum_{j\in S_{2}}\Pr\left\{  \sum_{i\in N_{j}}%
c_{ij}\left(  \xi_{i}-\xi_{i}^{\prime}\right)  =0\right\} \\
&  =\sum_{j\in S_{2}}\mbox{O}\left(  \frac{\sum_{i\in N_{j}}(1-\overline
{\kappa}(i))\overline{\kappa}_{\Delta}(i)}{\left\{  \sum_{i\in N_{j}}\left[
1-\overline{\kappa}_{\Delta}(i)\right]  \right\}  ^{3/2}}\right)  ,
\end{align*}
where $\overline{\kappa}(i)$ and $\overline{\kappa}_{\Delta}(i)$ are the jumps
associated with $\xi_{i}-\xi_{i}^{\prime}$. By Markov's inequality, we obtain
\[
\Pr\left\{  \zeta>\frac{|S_{2}|}{2}\right\}  \leq\frac{2}{|S_{2}|}%
\mathbb{E}\left(  \zeta\right)  =\frac{1}{|S_{2}|}\sum_{j\in S_{2}%
}\mbox{O}\left(  \frac{\sum_{i\in N_{j}}(1-\overline{\kappa}(i))\overline
{\kappa}_{\Delta}(i)}{\left\{  \sum_{i\in N_{j}}\left[  1-\overline{\kappa
}_{\Delta}(i)\right]  \right\}  ^{3/2}}\right)  .
\]

For $M:=\{j\in S_{2}:\eta_{j}\neq0\}$, we note that (i) $M$ is a random set
which depends only on $X,X^{\prime}$ and (ii) $|M|\geq|S_{2}|/2$ whenever
$\zeta\leq|S_{2}|/2$. Thus for a given realization $x$, $x^{\prime}$ of
$X,X^{\prime}$ respectively, we have
\[
\Pr\left\{  \varphi(x,Y)-\varphi(x^{\prime},Y)\in J\Bigm|\zeta\leq\frac
{|S_{2}|}{2}\right\}  =\Pr\left\{  2\sum_{j\in S_{2}}\xi_{j}\eta_{j}\in
J^{\prime}\Bigm|\zeta\leq\frac{|S_{2}|}{2}\right\}  ,
\]
where $J^{\prime}=[-g(x,x^{\prime})-\delta/2,-g(x,x^{\prime})+\delta/2]$. Then
by the Kolmogorov--Rogozin inequality,
\[
\Pr\left\{  \varphi(x,Y)-\varphi(x^{\prime},Y)\in J\Bigm|\zeta\leq\frac
{|S_{2}|}{2}\right\}  =\mbox{O}\left(  \frac{\sum_{j\in M(x,x^{\prime}%
)}(1-\kappa(j))\kappa_{\Delta}(j)}{\left\{  \sum_{j\in M(x,x^{\prime})}\left[
1-\kappa_{\Delta}(j)\right]  \right\}  ^{3/2}}\right)  ,
\]
where $M(x,x^{\prime})$ is the set $M$ obtained for the realization $x$,
$x^{\prime}$ of $X,X^{\prime}$. So
\begin{align*}
\MoveEqLeft[4]\Pr\left\{  \varphi(X,Y)-\varphi(X^{\prime},Y)\in J\Bigm|\zeta
\leq\frac{|S_{2}|}{2}\right\}  =\\
&  =\mathbb{E}\left(  \Pr\left\{  \varphi(X,Y)-\varphi(X^{\prime},Y)\in
J\Bigm|\zeta\leq\frac{|S_{2}|}{2},\;X,X^{\prime}\right\}  \right) \\
&  =\mathbb{E}\left(  \mbox{O}\left(  \sup_{D\subset S_{2},|D|\geq|S_{2}%
|/2}\frac{\sum_{j\in D}(1-\kappa(j))\kappa_{\Delta}(j)}{\left\{  \sum_{j\in
D}\left[  1-\kappa_{\Delta}(j)\right]  \right\}  ^{3/2}}\right)  \right)  \ \\
&  =\mbox{O}\left(  \sup_{D\subset S_{2},|D|\geq|S_{2}|/2}\frac{\sum_{j\in
D}(1-\kappa(j))\kappa_{\Delta}(j)}{\left\{  \sum_{j\in D}\left[
1-\kappa_{\Delta}(j)\right]  \right\}  ^{3/2}}\right)  .
\end{align*}
Hence
\[
\Pr\left\{  \varphi=x\right\}  =\mbox{O}\left(  \left[  \frac{1}{|S_{2}|}%
\sum_{j\in S_{2}}\left(  \frac{\sum_{i\in N_{j}}(1-\overline{\kappa
}(i))\overline{\kappa}_{\Delta}(i)}{\left\{  \sum_{i\in N_{j}}\left[
1-\overline{\kappa}_{\Delta}(i)\right]  \right\}  ^{3/2}}\right)
+\sup_{D\subset S_{2},|D|\geq|S_{2}|/2}\frac{\sum_{j\in D}(1-\kappa
(j))\kappa_{\Delta}(j)}{\left\{  \sum_{j\in D}\left[  1-\kappa_{\Delta
}(j)\right]  \right\}  ^{3/2}}\right]  ^{1/2}\right)  .
\]
\hfill$\blacksquare$

\begin{remark}
\label{RemWig01} a) If $\kappa_{\Delta}(i)<\kappa<1$ for all $i$,
$|S_{1}|=|S_{2}|=n/2$, and $|N_{j}|\geq n^{1-\varepsilon}$ for all $j$ and
$\varepsilon>0$,
\[
\Pr\left\{  \varphi=x\right\}  = \mbox{O} \left(  \left[  \frac{\kappa}%
{\sqrt{(1-\kappa)^{3} n^{1-\varepsilon}}}\right]  ^{1/2} \right)  .
\]

b) Lemma \ref{LemQuaCon} holds when $s$ of the random variables $\xi
_{1},\ldots,\xi_{n}$ \ are degenerate for some $1\leq s<n$, in this situation,
$n$ is replaced by $n-s$. The contribution to the bound of the concentration
inequality is only provided by the non-degenerate random variables.
\end{remark}

\section{Proofs in the Ginibre case}

We start with an extension of a result by Slinko \cite{Si01}, who treated the
case of a discrete uniform distribution with parameter $1/q$ with
$q\in\mathbb{Z}^{+}$. Throughout this section, all our random variables
satisfy
\[
\sup_{x\in\mathbb{R}}\Pr\{X=x\}\leq\kappa_{\Delta}(X)<\kappa<1.
\]

\begin{lemma}
\label{CorGin01} Let $k\leq m$ and let $A\in\mathbb{R}^{m\times k}$ be a
(deterministic) matrix with $\mbox{rank}(A)=k$. If $b\in\mathbb{R}^{m}$ is a
random vector whose entries are independent random variables, then
\[
\Pr\left\{  \mbox{rank}(A, b)=k\right\}  \leq\kappa^{m-k}.
\]

\end{lemma}

\noindent\textbf{Proof.} Since $\mbox{rank}(A)=k$, we can decompose $[A\; b]$
in the following way
\[
[A\; b]=\left(
\begin{array}
[c]{cc}%
A_{k} & b_{k}\\
A_{m-k} & b_{m-k}%
\end{array}
\right)  ,
\]
where $A_{k}\in\mathbb{R}^{k\times k}$, $A_{m-k}\in\mathbb{R}^{(m-k)\times k}%
$, $b_{k}\in\mathbb{R}^{k}$ and $b_{m-k}\in\mathbb{R}^{m-k}$. We note $A_{k}$
is an invertible matrix. We have that there exists a random matrix $\Delta
\in\mathbb{R}^{k}$ such that $A_{k}\Delta= b_{k}$ and $A_{m-k}\Delta= b_{m-k}%
$, then $A_{m-k}A^{-1}_{k} b_{k} =b_{m-k}$. So
\begin{align*}
\Pr\left\{  r(A, b) =k\right\}   &  \leq\Pr\left\{  A_{m-k}A^{-1}_{k} b_{k}
=b_{m-k} \right\} \\
&  = \mathbb{E}\left\{  \Pr\left\{  A_{m-k}A^{-1}_{k} b_{k} =b_{m-k} \left|
A_{m-k}A^{-1}_{k} b_{k}\right.  \right\}  \right\} \\
&  \leq\kappa^{m-k},
\end{align*}
the last line being due to the independence of every entry in $b_{m-k}$%
.\hfill$\blacksquare$

\begin{lemma}
\label{LemGin01} Let $k\leq m$ and let $A\in\mathbb{R}^{m\times k}$ be a
random matrix (whose entries are independent random variables). Then
\[
\Pr\left\{  \mbox{rank}(A)<k\right\}  <\frac{\kappa}{1-\kappa} \kappa^{m-k}.
\]

\end{lemma}

\noindent\textbf{Proof.} We note that if $A=[a_{1} |\cdots| a_{k}]$, $a_{i}%
\in\mathbb{R}^{m}$ $i=1,\ldots,k$, then
\begin{align*}
\Pr\{\mbox{rank}(A)=k\}  &  = \Pr\{a_{1}\notin\{0\}, a_{2}\notin%
\mbox{span}\{a_{1}\},\ldots, a_{k}\notin\mbox{span}\{a_{1},a_{2}%
,\ldots,a_{k-1}\} \}\\
&  = \Pr\{a_{1}\notin\{0\}\} \prod_{i=2}^{k} \Pr\{E_{i}\}
\end{align*}
where we use the notation $\mbox{span}\{\cdot\}$ for the space generated by
some vectors and
\[
E_{i}=\{a_{i}\notin\mbox{span}\{a_{1},a_{2},\ldots,a_{i-1}\} \left|
a_{1}\notin\{0\}, a_{2}\notin\mbox{span}\{a_{1}\},\ldots, a_{i-1}%
\notin\mbox{span}\{a_{1},a_{2},\ldots,a_{i-2}\}\right.  \}.
\]
Hence by Corollary \ref{CorGin01} and the Weierstrass product inequality,
\[
\Pr\{\mbox{rank}(A)=k\} \geq\prod_{i=0}^{k-1}(1-\kappa^{m-i}) \geq1-\sum
_{i=0}^{k-1} \kappa^{m-i} = 1 - \frac{\kappa}{1-\kappa} \kappa^{m-k}.
\]
\hfill$\blacksquare$

We consider the following concept used by Koml\'{o}s \cite{Bollobas1985}. Let
$S=\{v_{1},\ldots,v_{n}\}$ be a set of vectors. Let us define the strong rank
of $S$, denoted $sr(S)$, to be $n$ if $S$ is a set of linearly independent
vectors, and $k$ if any $k$ of the $v_{i}$s are linearly independent but some
$k+1$ of the vectors are linearly dependent. For a matrix $A$, we denote the
strong rank of the system of columns and the strong rank of the system of rows
by $sr_{c}(A)$ and $sr_{r}(A)$, respectively.

\begin{remark}
\label{RemGin02} (a) Let $A$ be an $m\times n$ random matrix with all entries
independent random variables. It follows immediately from Lemma~\ref{LemGin01}
that
\[
\Pr\left\{  sr_{c}(A)<k\right\}  \leq\binom{n}{k}\frac{\kappa}{1-\kappa}%
\kappa^{m-k}%
\]
(b) For every $\kappa$ and $0<\alpha\leq1$ there exists $\beta>0$ which
satisfies
\begin{equation}
\frac{h(\beta)}{\log_{2}\kappa}+\beta<\alpha<1, \label{LemGin02Eq01}%
\end{equation}
where $h(x)=-x\log_{2}(x)-(1-x)\log_{2}(1-x)$ is the entropy function. Indeed,
let
\[
g(x)=\frac{h(x)}{\log\kappa}+x.
\]
Now, since the function $g$ is continuous and $g(0)=0$, there exists a
positive number $\beta>0$ such that $g(\beta)<\alpha<1$.

c) We note from (a) and (b) that if $m=\lfloor\alpha n\rfloor$ and
$k=\lceil\beta n\rceil$, then
\[
\Pr\left\{  \mbox{rank}(A)<\lceil\beta n\rceil\right\}  <\binom{n}{\lceil\beta
n\rceil}\frac{\kappa}{1-\kappa}\kappa^{\lfloor\alpha n\rfloor-\lceil\beta
n\rceil}<\frac{\kappa}{1-\kappa}2^{n(h(\beta)-(\alpha-\beta)\log_{2}(\kappa
))}<\frac{\kappa}{1-\kappa}2^{-n\gamma_{\kappa}},
\]
where we use $\binom{n}{\beta n}<2^{nh(\beta)}$ and $\gamma_{\kappa}$ is a
positive constant which depends on $\kappa$.
\end{remark}

\begin{lemma}
\label{LemGin03} Let $v_{1},v_{2},\ldots,v_{k}\in\mathbb{R}^{m}$ be
(deterministic) linearly independent vectors. Let $B=[v_{1}| \ldots| v_{k} ]$
and $sc_{r}(B)=s$. Then for a random vector $a\in\mathbb{R}^{m}$ whose entries
are independent random variables,
\[
\Pr\left\{  \mbox{rank}(v_{1},v_{2},\ldots,v_{k},a)=k\right\}  <C_{1}
\kappa^{m-k}s^{-1/2}.
\]

\end{lemma}

\noindent\textbf{Proof.} Although simple, for the sake of completeness we
include the proof. Let $b_{1},b_{2},\ldots,b_{m}$ be the rows of $B$. Without
loss of generality we may assume that $b_{1},b_{2},\ldots,b_{k}$ are linearly
independent and that all other rows are linear combination of them. We have
\[
\sum_{i=1}^{k}\beta_{i}^{(r)}b_{i}=b^{(r)}%
\]
for $r=k+1,\ldots,m$. As $sc_{r}(B)=s$, at least $s$ of the coefficients
$\beta_{1}^{(r)}\ldots,\beta_{k}^{(r)}$ are nonzero.

Now, since we consider the event $[\mbox{rank}(v_{1},v_{2},\ldots,v_{k}%
,a)=k]$, we have
\[
\sum_{j=1}^{k}\alpha_{j}v_{j}=a
\]
for some $\alpha_{1}\ldots,\alpha_{k}$ not all zero. In particular $\sum
_{j=1}^{k}\alpha_{j}v_{k+1,j}=a_{k+1}$, where $a_{k+1}$ is the $(k+1)$th entry
of $a$. But
\[
a_{k+1}=\sum_{j=1}^{k}\alpha_{j}v_{k+1,j}=\sum_{j=1}^{k}\alpha_{j}\left(
\sum_{i=1}^{k}\beta_{i}^{(k+1)}v_{i,j}\right)  =\sum_{i=1}^{k}\beta
_{i}^{(k+1)}\left(  \sum_{j=1}^{k}\alpha_{j}v_{i,j}\right)  =\sum_{i=1}%
^{k}\beta_{i}^{(k+1)}a_{i}.
\]
From the above and the independence of the entries of $a$,
\begin{align*}
\Pr\left\{  \mbox{rank}(v_{1},v_{2},\ldots,v_{k},a)=k\right\}   &  \leq
\Pr\left\{  \sum_{i=1}^{k}\beta_{i}^{(r)}a_{i}=a_{r},r=k+1,\ldots,m\right\} \\
&  =\mathbb{E}\left\{  \Pr\left\{  \sum_{i=1}^{k}\beta_{i}^{(r)}a_{i}%
=a_{r},r=k+1,\ldots,m\left\vert {}\right.  a_{1},\ldots,a_{k}\right\}
\right\} \\
&  =\mathbb{E}\left\{  \Pr\left\{  \sum_{i=1}^{k}\beta_{i}^{(m)}a_{i}%
=a_{m}\left\vert {}\right.  a_{1},\ldots,a_{k}\right\}  \prod_{l=k+1}^{m-1}%
\Pr\left\{  \sum_{i=1}^{k}\beta_{i}^{(l)}a_{i}=a_{l}\left\vert {}\right.
a_{1},\ldots,a_{k}\right\}  \right\} \\
&  \leq\mathbb{E}\left\{  \kappa^{m-k-1}\Pr\left\{  \sum_{i=1}^{k}\beta
_{i}^{(m)}a_{i}=a_{m}\left\vert {}\right.  a_{1},\ldots,a_{k}\right\}
\right\} \\
&  =\kappa^{m-k-1}\Pr\left\{  \sum_{i=1}^{k}\beta_{i}^{(m)}a_{i}=a_{m}\right\}
\\
&  \leq C_{1}\kappa^{m-k}s^{-1/2},
\end{align*}
the last line being due to Lemma \ref{LemLinCon} and Remark \ref{RemLinConIne}%
.\hfill$\blacksquare$

\textbf{Proof of Theorem \ref{Univer}.a}. Let $\alpha\in(0,1)$ and $\beta>0$
be as in Equation (\ref{LemGin02Eq01}) and let $n_{0}=\lfloor\alpha n\rfloor$.
Let $B$ be the $n_{0}\times n$ matrix whose columns are the first $n_{0}$
columns of $G_{n}$.

From Lemma \ref{LemGin01} we can assume that $B$ has full rank. Since
\[
\Pr\{\mbox{rank}(G_{n})=n\} = \Pr\{\mbox{rank}(G_{n})=n, sr_{r}(B)<\beta n\} +
\Pr\{\mbox{rank}(G_{n})=n, sr_{r}(B)\geq\beta n\},
\]
by Lemma \ref{LemGin03} and Remark \ref{RemGin02}, we have%

\[
\Pr\{\mbox{rank}(G_{n})=n\} \geq\prod_{i=1}^{n-n_{0}} \left(  1- C_{1}(\beta
n)^{-1/2}\kappa^{i}\right)  \geq1 - \frac{C_{1}}{1-\kappa} (\beta n)^{-1/2},
\]
which proves Theorem \ref{Univer}.a. \hfill$\blacksquare$

\textbf{Proof of Proposition \ref{ProGin01}}. Let $F_{1}$ be a distribution
function whose biggest jump is $\kappa_{1}$. We take $m_{n}=1$ and $\delta
_{1}=\kappa_{1}/2$, then $\Pr\{ G_{m_{1}} \mbox{ has full rank }\}>1-\delta
_{1}$. Now, let $F_{n}$ be a distribution function whose biggest jump is
$\kappa_{n}$. By Lemma 2 in \cite{Komlos1968}, there is $m_{n}\geq m_{n-1}$
and $\delta_{n}\leq1/n\leq$ for $n>1$ such that
\[
\Pr\{ G_{m_{n}} \mbox{ has full rank }\}>1-\delta_{n}%
\]
where the entries of $G_{m_{n}}$ have the same distribution and $\delta_{n}%
\to0$ as $n\to\infty.$ \hfill$\blacksquare$

\section{Proofs in the Wigner case}

\label{Section:Wigner}

Following the terminology introduced in Costello, Tao and Vu
\cite{CostelloTaoVu2006}$,$ given $n$ vectors $\left\{  v_{1},\ldots
,v_{n}\right\}  $, a linear combination of the $v_{i}$s is a vector
$v=\sum_{i=1}^{n}c_{i}v_{i}$, where the $c_{i}$ are real numbers. We say that
a linear combination vanishes if $v$ is the zero vector. A vanishing linear
combination has degree $k$ if exactly $k$ among the $c_{i}$ are nonzero.

A singular $n\times n$ matrix is called \textit{normal} if its row vectors do
not admit a non-trivial vanishing linear combination with degree less than
$n^{1-\varepsilon}$ for a given $\varepsilon\in(0,1)$. Otherwise it is said
that the matrix is \textit{abnormal}. Furthermore, a row of an $n\times n$
non-singular matrix is called \textit{good} if its exclusion leads to an
$(n-1)\times n$ matrix whose column vectors admit a non-trivial vanishing
linear combination with degree at least $n^{1-\varepsilon}$ (in fact, there is
exactly one such combination, up to scaling, as the rank of this $(n-1)\times
n$ matrix is $n-1$). A row is said to be \textit{bad} otherwise. Finally, an
$n\times n$ non-singular matrix $A$ is \textit{perfect} if every row in $A$ is
good. If a non-singular matrix is not perfect, it is called \textit{imperfect}.

For the proof of Theorem~\ref{Univer}.b, we first present three lemmas which
generalize results in \cite{CostelloTaoVu2006} for Wigner matrices $W_{n}%
=(\xi_{ij})$ with independent entries which need not be identically
distributed and the appropriate estimates in these new cases are found in
terms of the size of the biggest jump of the distribution functions governing
the entries under the hypothesis $\kappa_{\Delta}(i)<\kappa<1$. We also obtain
a better rate of convergence, which is universal. The proofs we give follow
ideas in \cite{CostelloTaoVu2006} but also take into account the size of the
biggest jump.

\begin{lemma}
\label{LemWig01} Let $\varepsilon\in(0,1)$, then for all $n$ large
\begin{equation}
\Pr\left\{  \mbox{$W_n$ is singular and abnormal}\right\}  \leq\kappa
^{(n-n^{1-\varepsilon})/2} \label{LemWig11}%
\end{equation}
and
\begin{equation}
\Pr\left\{  \mbox{$W_n$ is non-singular and imperfect}\right\}  \leq
\kappa^{(n-n^{1-\varepsilon})/2}. \label{LemWig12}%
\end{equation}

\end{lemma}

\textbf{Proof.} If $W_{n}$ is singular and abnormal the row vectors of $W_{n}$
admit a non-trivial vanishing linear combination with degree at most
$N:=n^{1-\varepsilon}$. For $i=1,\ldots,N$, we have that if $i=1$, there is a
row of $W_{n}$ that contains only zeros, and if $i>1$, the $i$th row is a
linear combination of the first $i-1$ rows of $W_{n}$ that are linearly
independent. We denote by $D(n,i)$ this last event and by $T_{i-1}$ the upper
triangular part of $W_{n}$ until the row $i-1$ (included). The linear
dependence of the $i$th row of $W_{n}$ with the $i-1$ rows of $W_{n}$ is
determined only by its last $n-i+1$ entries. Then by the stochastic
independence of $T_{i-1}$ with the last $n-i+1$ entries of the row $i$
\begin{align*}
\Pr\left\{  \mbox{$W_n$ is singular and abnormal}\right\}   &  \leq\sum
_{i=1}^{N}\binom{n}{i}\Pr\left\{  D(n,i)\right\}  \leq\sum_{i=1}^{N}\binom
{n}{i}\mathbb{E}\left\{  \Pr\left\{  \left.  D(n,i)\right\vert T_{i-1}%
\right\}  \right\} \\
&  \leq\sum_{i=1}^{N}n^{N}\kappa^{n-N+1}=Nn^{N}\kappa^{n-N+1},
\end{align*}
and for all $n$ large,
\[
\Pr\left\{  \mbox{$W_n$ is singular and abnormal}\right\}  \leq\kappa
^{\frac{3}{4}(n-n^{1-\varepsilon})}\leq\kappa^{\frac{1}{2}(n-n^{1-\varepsilon
})}.
\]

Now, we consider the case when $W_{n}$ is non-singular and imperfect. We can
suppose that the last row of $W_{n}$ is the bad row. The $(n-1)\times
n$-matrix obtained has rank $n-1$, hence there is a unique column that admits
a non-trivial vanishing linear combination with degree at most
$n^{1-\varepsilon}$. Then the last $n-k-1$ entries of this column are
completely determined by its $k$ first entries and $k$ linearly independent
columns, for $1\leq k\leq n^{1-\varepsilon}$. Since we can choose this bad
row, we have as above for $n$ large
\[
\Pr\left\{  \mbox{$W_n$ is non-singular and imperfect}\right\}  \leq
n\kappa^{\frac{3}{4}(n-1-(n-1)^{1-\varepsilon})}\leq\kappa^{\frac{1}%
{2}(n-n^{1-\varepsilon})}.
\]
\hfill$\blacksquare$

\begin{lemma}
\label{LemWig02} Let $A$ be a deterministic $n\times n$ singular normal
matrix. Then
\[
\Pr\left\{  \mbox{rank}(W_{n+1})-\mbox{rank}(W_{n})<2\left|  W_{n}=A\right.
\right\}  = \mbox{O}_{\varepsilon}\left(  \frac{\kappa}{\sqrt{n^{1-\varepsilon
}(1-\kappa)^{3}}}\right)  .
\]

\end{lemma}

\textbf{Proof.} Since $r:=\mbox{rank}(A)<n$, without loss of generality it is
possible to suppose that the first $r$ rows of $A$ are linearly independent.
If $v_{1},\ldots,v_{r}$ are the first rows of $A$, then $v_{n}=\sum_{i=1}%
^{r}\alpha_{i}v_{i}$, and as $A$ is normal, the number of coefficients in this
linear combination is at least $n^{1-\varepsilon}$. If it does not hold that
$\xi_{n}=\sum_{i=1}^{r}\alpha_{i}\xi_{i}$, where $\xi_{i}$ are entries of the
last column of $W_{n+1}$, by symmetry of $W_{n+1}$ we have
$\mbox{rank}(W_{n+1})=\mbox{rank}(A)+2$. Hence
\begin{align*}
\Pr\left\{  \mbox{rank}(W_{n+1})-\mbox{rank}(W_{n})<2\left\vert W_{n}
=A\right.  \right\}   &  \leq\Pr\left\{  \xi_{n}=\sum_{i=1}^{r}\alpha_{i}%
\xi_{i}\right\} \\
&  = \mbox{O}_{\varepsilon}\left(  \frac{\kappa}{\sqrt{n^{1-\varepsilon
}(1-\kappa)^{3}}}\right)  .
\end{align*}
The last expression follows from Lemma~\ref{LemLinCon}.\hfill$\blacksquare$

\begin{lemma}
\label{LemWig03} Let $A$ be a deterministic $n\times n$ non-singular perfect
symmetric matrix. Then
\[
\Pr\left\{  \mbox{rank}(W_{n+1})=n\left\vert W_{n}=A\right.  \right\}
=\mbox{O}_{\varepsilon}\left(  \left[  \frac{\kappa}{\sqrt{n^{1-\varepsilon
}(1-\kappa)^{3}}}\right]  ^{1/2}\right)  .
\]

\end{lemma}

\textbf{Proof.} If $\mbox{rank}(W_{n+1})=n$, then $\det(W_{n+1})=0$, and we
have
\[
0=\det(W_{n+1})=(\det A)\xi_{n+1}+\sum_{i=1}^{n}\sum_{j=1}^{n}c_{ij}\xi_{i}%
\xi_{j},
\]
where $\xi_{i}$ are entries of the last column of $W_{n+1}$ and its transpose,
and the $c_{ij}$ are the cofactors of $A$. Since $A$ is perfect, when we
eliminate the $i$th row of $A$, the columns of the matrix thus obtained admit
a vanishing linear combination of degree at least $n^{1-\varepsilon}$. When
the column $j$ is selected, where $j$ is the index of a non-zero coefficient
in this linear combination, we obtain an $(n-1)\times(n-1)$ non-singular
matrix since there are at least $n^{1-\varepsilon}$ indices $i$ such that
there are at least $n^{1-\varepsilon}$ indices $j$ with $c_{i,j}\neq0$. Taking
the partition of $\{1,2,\ldots,n\}$ as $S_{1}=\{1,2,\ldots,\lfloor
n/2\rfloor\}$ and $S_{2}=\{1,2,\ldots,n\}-S_{1}$, by Remark \ref{RemWig01}
\begin{align*}
\Pr\left\{  \mbox{rank}(W_{n+1})=n\left\vert W_{n}=A\right.  \right\}   &
\leq\Pr\left\{  (\det A)\xi_{n+1}+\sum_{i=1}^{n}\sum_{j=1}^{n}c_{ij}\xi_{i}%
\xi_{j}=0\right\} \\
&  =\mathbb{E}\left(  \Pr\left\{  (\det A)\xi_{n+1}+\sum_{i=1}^{n}\sum
_{j=1}^{n}c_{ij}\xi_{i}\xi_{j}=0\left\vert {}\right.  \xi_{n+1}\right\}
\right) \\
&  =\mathbb{E}\left(  \mbox{O}_{\varepsilon}\left(  \left[  \frac{\kappa
}{\sqrt{n^{1-\varepsilon}(1-\kappa)^{3}}}\right]  ^{1/2}\right)  \right) \\
&  =\mbox{O}_{\varepsilon}\left(  \left[  \frac{\kappa}{\sqrt{n^{1-\varepsilon
}(1-\kappa)^{3}}}\right]  ^{1/2}\right)  .
\end{align*}
\hfill$\blacksquare$

Now we consider the discrete stochastic process
\[
X_{n}=\left\{
\begin{array}
[c]{ll}%
0 & \mbox{if rank$(W_n)=n$}\\
\left(  \kappa^{-1/8}\right)  ^{n-\mbox{rank}(W_{n})} & \mbox{if
rank$(W_n)<n$},
\end{array}
\right.
\]
for which we can prove the following result.

\begin{proposition}
\label{ProWig01}
\[
\mathbb{E}\left(  X_{n}\right)  =\mbox{O}_{\varepsilon}\left(  \left[
\frac{\kappa}{\sqrt{n^{1-\varepsilon}(1-\kappa)^{3}}}\right]  ^{1/2}\right)
.
\]

\end{proposition}

\textbf{Proof.} For $j=0,\ldots,n$, write $A_{j}=\{\mbox{rank}(W_{n})=n-j\}$
and let $1+\gamma=\kappa^{-1/8}$. We have
\begin{align*}
\mathbb{E}\left(  X_{n}\right)   &  =\sum_{j=1}^{n}(1+\gamma)^{j}\Pr\left\{
A_{j}\right\} \\
&  =\sum_{j=1}^{n}(1+\gamma)^{j}\Pr\left\{  A_{j}%
\mbox{, $W_n$ normal}\right\}  +S_{1},
\end{align*}
where
\[
S_{1}=\sum_{j=1}^{n}(1+\gamma)^{j}\Pr\left\{  A_{j}%
\mbox{, $W_n$ abnormal}\right\}  .
\]
By Lemma \ref{LemWig01},
\begin{align*}
S_{1}  &  \leq\sum_{j=1}^{n}(1+\gamma)^{j}\kappa^{(n-n^{1-\varepsilon)}/2}\\
&  \leq\kappa^{(n-n^{1-\varepsilon)}/2}\sum_{j=1}^{n}(1+\gamma)^{j}\\
&  \leq\frac{1-(\kappa^{-1/8})^{n+1}}{1-\kappa^{-1/8}}\kappa
^{(n-n^{1-\varepsilon)}/2}\\
&  =C\kappa^{(3n-4n^{1-\varepsilon})/8}%
\end{align*}
for some constant $C>0$.

So
\begin{equation}
\label{eq1}\mathbb{E}\left(  X_{n}\right)  = \sum_{j=1}^{n} (1+\gamma)^{j}
\Pr\left\{  A_{j}\mbox{, $W_n$ normal} \right\}  + \mbox{O}_{\varepsilon
}\left(  \kappa^{(3n-4n^{1-\varepsilon})/8}\right)  .
\end{equation}

On the other hand,
\[
\mathbb{E}\left(  X_{n+1}\right)  = S_{2}+S_{3}+S_{4}+S_{5},
\]
where
\begin{align*}
S_{2}  &  = \mathbb{E}\left(  X_{n+1}\left\vert A_{0}%
,\mbox{$W_n$ perfect}\right.  \right)  \Pr\left\{  A_{0}%
\mbox{, $W_n$ perfect}\right\} \\
S_{3}  &  = \mathbb{E}\left(  X_{n+1}\left\vert A_{0}%
,\mbox{$W_n$ imperfect}\right.  \right)  \Pr\left\{  A_{0}%
\mbox{, $W_n$ imperfect}\right\} \\
S_{4}  &  = \sum_{j=1}^{n}\mathbb{E}\left(  X_{n+1}\left\vert A_{j}%
,\mbox{$W_n$
normal}\right.  \right)  \Pr\left\{  A_{j}\mbox{, $W_n$ normal}\right\} \\
S_{5}  &  = \sum_{j=1}^{n}\mathbb{E}\left(  X_{n+1}\left\vert A_{j}%
,\mbox{$W_n$
abnormal}\right.  \right)  \Pr\left\{  A_{j}\mbox{, $W_n$ abnormal}\right\}  .
\end{align*}

By Lemma \ref{LemWig03} and the fact that rank$(W_{n})=n$,
\begin{align*}
S_{2}  &  \leq(\kappa^{-1/8})^{n+1-n}\Pr\{\mbox{rank}(W_{n+1})=n\left\vert
W_{n} \mbox{ is perfect and non-singular}\right.  \}\\
&  = \mbox{O}_{\varepsilon}\left(  \left[  \frac{\kappa}{\sqrt
{n^{1-\varepsilon}(1-\kappa)^{3}}}\right]  ^{1/2}\right)  .
\end{align*}

On the other hand, Lemma \ref{LemWig01} and the definition of $X_{n+1}$ give
\[
S_{3}\leq(\kappa^{-1/8})^{n+1}\kappa^{(n-n^{1-\varepsilon})/2}%
=\mbox{O}_{\varepsilon}\left(  \kappa^{(3n-4n^{1-\varepsilon})/8}\right)  .
\]

Using again Lemma \ref{LemWig01} and the definition of $A_{j}$,
\[
S_{5}\leq\sum_{j=1}^{n}(\kappa^{-1/8})^{j+1}\kappa^{(n-n^{1-\varepsilon}%
)/2}=\mbox{O}_{\varepsilon}\left(  \kappa^{(3n-4n^{1-\varepsilon})/8}\right)
.
\]

If rank$(W_{n})=n-j$, then rank$(W_{n+1})$ is equal to $n-j+2$ or $n-j$ since
$W_{n+1}$ is a symmetric matrix. By Lemma \ref{LemWig02} and for $n$
sufficiently large,
\begin{align*}
\mathbb{E}\left(  X_{n+1}\left\vert A_{j},\mbox{$W_n$ normal}\right.  \right)
&  =(1+\gamma)^{j+1}\Pr\{\mbox{rank}(W_{n+1})=\mbox{rank}(W_{n})\left\vert
\mbox{$W_n$ normal and singular}\right.  \}\\
&  +(1+\gamma)^{j-1}\\
&  =(1+\gamma)^{j}\left(  (1+\gamma)^{-1}+\mbox{O}_{\varepsilon}\left(
\frac{\kappa}{\sqrt{n^{1-\varepsilon}(1-\kappa)^{3}}}\right)  \right) \\
&  \leq\alpha(1+\gamma)^{j}%
\end{align*}
for some $\alpha<1$.

Then we have
\[
\mathbb{E}(X_{n+1})=\alpha\sum_{j=1}^{n}(1+\gamma)^{j}\Pr\{A_{j}%
,W_{n}\mbox{ normal}\}+\mbox{O}_{\varepsilon}\left(  f(\kappa,n)\right)  ,
\]
where
\[
f(\kappa,n):=\frac{\kappa^{\frac{3}{8}n-\frac{1}{2}n^{1-\varepsilon}}}%
{\kappa(1-\kappa)}+\left[  \frac{\kappa}{\sqrt{n^{1-\varepsilon}(1-\kappa
)^{3}}}\right]  ^{1/2}.
\]

Using (\ref{eq1})
\[
\mathbb{E}( X_{n+1}) \leq\alpha\mathbb{E}( X_{n}) + \mbox{O}_{\varepsilon
}\left(  f(\kappa,n)\right)  ,
\]
so
\[
\mathbb{E}( X_{n+1}) \leq\alpha^{n} \mathbb{E}( X_{1}) + \mbox{O}_{\varepsilon
}\left(  f(\kappa,n)\right)  .
\]

This proves the proposition.\hfill$\blacksquare$

\textbf{Proof of Theorem \ref{Univer}.b}. By Markov's inequality,
\begin{align}
\Pr\left\{  \mbox{rank}(W_{n})<n\right\}   &  =\Pr\left\{  X_{n}\geq1\right\}
\nonumber\\
&  \leq\mathbb{E}\left(  X_{n}\right) \nonumber\\
&  =\mbox{O}_{\varepsilon}\left(  \left[  \frac{\kappa}{\sqrt{n^{1-\varepsilon
}(1-\kappa)^{3}}}\right]  ^{1/2} \right)  , \label{Wiglast}%
\end{align}
where we have used Proposition \ref{ProWig01}. \hfill$\blacksquare$

\section*{Acknowledgments}

The authors would like to thank the constructive and useful suggestions
provided by the referees, AE and Editor, which improved the manuscript. Rahul
Roy wants to thank CIMAT for the warm hospitality he received during his
visits. The work of Paulo Manrique was supported by the Ph.D. Conacyt grant 210223.

\end{document}